\input amstex
\magnification=\magstep0
\documentstyle{amsppt}
\pagewidth{6.0in}
\input amstex
\topmatter
\title
Theorems, Problems and Conjectures
\endtitle
\author
Tewodros Amdeberhan\\
version 25 August, 2022
\endauthor
\affil
Department of Mathematics\\
Tulane University\\
New Orleans, LA 70118\\
tamdeber\@tulane.edu
\endaffil
\abstract
These notes are designed to offer some (perhaps new) codicils to related work, a list of problems and
conjectures seeking (preferably) combinatorial proofs. The main items are Eulerian polynomials
and hook/contents of Young diagram, mostly on the latter. We also have items on Frobenius theorem
and multi-core partitions; most recently, some problems on (what we call) colored overpartitions. Formulas analogues to or in the spirit of works by Han, Nekrasov-Okounkov and Stanley are distributed throughout. Concluding remarks are provided at the end in hopes of directing the interested researcher, properly. The newly added problem is on chromatic polynomials.
\endabstract
\endtopmatter
\def\({\left(}
\def\){\right)}

\document
\noindent
\bf Notations \rm
\bigskip

$p(n)=$ number of unrestricted partitions of $n$

$S_n=$ the symmetric group of permutation on $\{1,2,\dots,n\}$

$Inv(S_n)=\{\pi\in S_n: \pi^2=1\}=$ involutions


$S(n,k)=$ Stirling numbers of the second kind

$\widehat{B}_n(z):=\sum_{j=0}^nS(n,j)z^j=$ Bell polynomials

$\lambda\vdash n$ means $\lambda$ is a unrestricted partition of $n$; $\lambda^t=$ conjugate partition to $\lambda$

$L(\lambda)=$ length (number of parts) of a partition $\lambda$

$\Cal{H}(\lambda)=$ all hooks of cells in $\lambda$

$c_k(\pi)=$ number of cycles of $\pi$ of length $k$; $c_1(\pi)=tr(\pi)=$ trace of $\pi$

$\kappa(\pi)=$ number of cycles of a permutation $\pi$

$\Cal{O}(\pi)=$ number of odd-length cycles of a permutation $\pi$

$\Cal{E}(\pi)=$ number of even-length cycles of a permutation $\pi$

$l_k(\pi)=$ number of cycles, in $\pi\in S_n$, of length $k$

$c_u=$ the content of the cell $u$ in a partition $\lambda$

$h_u=h(i,j)$ the hook length of cell $u=(i,j)$ in a partition $\lambda$

$a_u=$ arm of a hook; $\ell_u=$ leg of a hook (hence $h_u=a_u+\ell_u+1$)

$f_{\lambda}=\dim\lambda=$ number of Standard Young Tableau of shape $\lambda$

$s_{\lambda}=$ Schur function

$e_{\lambda}=$ elementary symmetric functions

Partitions $\lambda$ will also be written in powers $\lambda=(1^{m_1}2^{m_2}\cdots)$

$(a;q)_n:=(1-a)(1-aq)(1-aq^2)\cdots(1-aq^{n-1})$ with $(a;q)_0:=1$

$(a;q)_{\infty}=(1-a)(1-aq)(1-aq^2)\cdots$; $\pmb{\binom{n}k}=\frac{(q;q)_n}{(q;q)_k(q;q)_{n-k}}$, the Gaussian polynomials.

\pagebreak
\bigskip
\noindent
\bf 1 EULERIANS \rm
\smallskip
\noindent
The \it Eulerian polynomials of type A \rm and type B \rm are defined, respectively, by
$$\sum_{k\geq0}(k+1)^nt^k=\frac{A_n(t)}{(1-t)^{n+1}}, \qquad\text{and} \qquad
\sum_{k\geq0}(2k+1)^nt^k=\frac{B_n(t)}{(1-t)^{n+1}}.$$
The following three results are easy to establish.
\bigskip
\noindent
\bf Lemma 1.1 \it If $A_n(t)=\sum_{j\geq0}A_{n,j}t^j$ and $B_n(t)=\sum_{j\geq0}B_{n,j}t^j$, then
the Eulerian numbers are linked by
$$2^nA_{n,k}=\sum_{j\geq0}\binom{n+1}{2k+1-j}A_{n,j} \qquad\text{and} \qquad
B_{n,k}=\sum_{j\geq0}\binom{n+1}{2k-j}A_{n,j}.$$ \rm
\bf Lemma 1.2 \it The following relation intertwines type A and type B
$$\sum_{k\geq0}\binom{a+b}kB_{k,b}=\sum_{k\geq0}\binom{a+b}k2^kA_{k,a-1}.$$ \rm
\bf Lemma 1.3 \it The following relates Stirling numbers of second kind with Eulerian polynomial by
$$A_k(t)=\sum_{j=1}^kj!S(k,j)t^{j-1}(1-t)^{k-j}.$$ \rm
\bf Definition 1.4 \rm A \it $q$-exponential function \rm $e(z;q)$ is defined as
$$e(z;q)=\sum_{n\geq0}\frac{z^n}{(q;q)_n}.$$
\bf Definition 1.5 \rm A \it $q$-Eulerian polynomial of type B \rm is defined by the generating function
$$\sum_{n\geq1}B_n(t,q)\frac{z^n}{(q;q)_n}=\frac{(e(z;q)-e(tz;q))(e(tz;q)+te(z;q))}{e(2tz;q)-te(2z;q)}.$$
\bf Definition 1.6 \rm The \it $q$-Eulerian numbers of type B \rm are defined by
$$B_n(t,q)=\sum_{k=0}^nB_{n,k}(q)t^k.$$
The next two statements reveal interesting properties of the $q$-Eulerians of type B.
\bigskip
\noindent
\bf Lemma 1.7 \it We have the symmetry
$$B_{n,k}(q)=B_{n,n-k}(q).$$ \rm
\bf Proof. \rm Substituting $t\rightarrow t^{-1}$ and $z\rightarrow tz$ implies $t^nB_n(t^{-1},q)=B_n(t,q)$.
The claim follows. $\square$
\bigskip
\noindent
\bf Conjecture 1.8 \it Let $a, b\in\Bbb{N}$ and denote $\alpha=a+b+1$. We have the symmetric relation
$$\binom{\alpha}a_q+\sum_k\binom{\alpha}k_q2^{\alpha-k}B_{k,b}(q)=
\binom{\alpha}b_q +\sum_k\binom{\alpha}k_q2^{\alpha-k}B_{k,a}(q).$$ \rm

\pagebreak
\bigskip
\noindent
\bf 2 HOOKS, ARMS and LEGS
\bigskip
\noindent
\bf Conjecture 2.1 \it
For a partition $\lambda\vdash n$ and its Young diagram, let $h_u, a_u, \ell_u$ denote the hook-length,
arm and legs of a cell $u\in\lambda$. Note $h_u=a_u+\ell_u+1$. Let $k_j(\lambda)=\#\{i: \lambda_i=j\}$.
Then, we have variants linked to Nekrasov-Okounkov hook formula \rm [8]
$$\sum_{\lambda\vdash n}\prod\Sb u\in\lambda\\ a_u=0\endSb\frac{h_u+t}{h_u}=
\sum_{\lambda\vdash n}\prod_{j\geq1}\binom{k_j+t}{k_j}=
\sum_{\lambda\vdash n}\prod\Sb u\in\lambda\endSb\frac{h_u^2+t}{h_u^2}=
\sum_{\lambda\vdash n}\prod\Sb u\in\lambda\\ \ell_u=0\endSb\frac{h_u+t}{h_u}.$$ \rm
\bf Conjecture 2.2(a) \it Each of the following polynomials, in $t$,  has only simple negative real roots
$$P_n(t)=\sum_{\lambda\vdash n}\prod\Sb u\in\lambda\endSb\frac{h_u^2+t}{h_u^2}.$$
\bf Conjecture 2.2(b) \it Each of the polynomials $P_n(t)$ is unimodal. \rm
\bigskip
\noindent
\bf  Problem 2.2 \it Denote $d_1(\lambda)=\#\{u\in\lambda: h_u=1\}$ and $b_n=\max\{d_1(\lambda): \lambda\vdash n\}$. Then, we have
$$\sum_{n\geq0}b_nx^n=\frac{x}{1-x}\frac{(x^2;x^2)_{\infty}^2}{(x;x)_{\infty}}.$$ \rm
\bf Lemma 2.3 \it If $\lambda\vdash n$, let $f_1(\lambda)=\#\{i:\lambda_i=1\}$ as the number of times
$1$ appears as a part of $\lambda$, $g_1(\lambda)$ be the number of distinct parts of $\lambda$,
$d_1(\lambda)=\#\{u\in\lambda: h_u=1\}$ as the number of times $1$ appears as a hook-length of
$\lambda$ and $p(k)=\#\{\lambda:\lambda\vdash k\}$ denotes the number of unrestricted partitions of $k$ .
Then, the following are equal:
$$\sum_{\lambda\vdash n}f_1(\lambda)=\sum_{\lambda\vdash n}g_1(\lambda)=\sum_{\lambda\vdash n}d_1(\lambda)
=\sum_{k=0}^np(k).$$
\bf Remark. \rm Parts of Lemma 2.3 appear in Stanley's EC-1, Exercise 1.26 (1999 first ed.)
or Exercise 1.80 (2nd edition).

\bigskip
\noindent
\bf 3 UNIFYING the KNOWNS \rm
\bigskip
\noindent
The main purpose in this section is to exploit one single formula (next lemma) so as we
connect certain known identities that enjoy hooks and contents representations. We emphasize that
the hook/content sums are not claimed to have been proven by this method; those are borrowed. \rm
\smallskip
\noindent
\bf Lemma 3.1 \it We have the generating function
$$(1-x)^{-t}(1-x^2)^{-\left(\frac{q-t}2\right)}=\sum_{n=0}^{\infty}\frac{x^n}{n!}
\sum_{\pi\in S_n}t^{\Cal{O}(\pi)}q^{\Cal{E}(\pi)}.\tag1$$ \rm
\bf Proof. \rm Standard exponential generating function techniques (see e.g. [11, eqn. (5.30)]). $\square$
\smallskip
\noindent
\bf Example 3.2 \rm Replace $q=t=tv$ in (1). The result relates to Stanley's "content Nekrasov-Okounkov" formula [10, Thm. 2.2].
$$\sum_{n=0}^{\infty}x^n\sum_{\lambda\vdash n}\prod_{u\in\lambda}\frac{(t+c_u)(v+c_u)}{h_u^2}=
\frac1{(1-x)^{tv}}=\sum_n\frac{x^n}{n!}\sum_{\pi\in S_n}(tv)^{\kappa(\pi)}.$$
\bf Example 3.3 \rm Replace $q=t^2$ in (1). \it Is this generating function \bf new\rm? \rm
$$\sum_nx^n\sum_{\lambda\vdash n}\prod_{u\in\lambda}\frac{t+c_u}{h_u}=
\frac1{(1-x)^t(1-x^2)^{\binom{t}2}}=\sum_n\frac{x^n}{n!}\sum_{\pi\in S_n}t^{\Cal{O}(\pi)+2\Cal{E}(\pi)}.$$ \rm
\bf Example 3.4 \rm Replace $x=yz, v=\frac1z$ in Example 3.2, and take the limit $z\rightarrow 0$.
$$\sum_ny^n\sum_{\lambda\vdash n}\prod_{u\in\lambda}\frac{t+c_u}{h_u^2}=e^{ty}=
\sum_n\frac{y^n}{n!}\sum\Sb\pi\in S_n\\  \pi=1\endSb t^{tr(\pi)}.$$
\bf Example 3.5 \rm Replace $t=ab, q=b$. Then $t=a, q=1$ and raise to the $b$ power to get left and right-sides.
$$\sum_n\frac{x^n}{n!}\sum_{\pi\in S_n}a^{\Cal{O}(\pi)}b^{\kappa(\pi)}=
\left(\sum_n\frac{x^n}{n!}\sum_{\pi\in S_n}a^{\Cal{O}(\pi)}\right)^b.$$
\bf Example 3.6 \rm  Put $x=yz, t=\frac1z, q=\frac{a^2}{z^2}$ in equation (1) to get
$$(1-yz)^{-\frac1z}(1-y^2z^2)^{-\left(\frac{a^2-z}{2z^2}\right)}=
\sum_{n=0}^{\infty}\frac{y^n}{n!}\sum_{\pi\in S_n}a^{2\Cal{E}(\pi)}z^{n-\Cal{O}(\pi)-2\Cal{E}(\pi)}.$$
Now, take the limit $z\rightarrow 0$ on both sides, and observe that $n=\Cal{O}(\pi)+2\Cal{E}(\pi)$ iff $\pi^2=1$, in which case $\Cal{E}(\pi)=c_2(\pi)$. The result relates to Han's formula [3]
$$\sum_{n=0}^{\infty}y^n\sum_{\lambda\vdash n}\prod_{u\in\lambda}
\frac{(1+a)^{h_u}+(1-a)^{h_u}}{(1+a)^{h_u}-(1-a)^{h_u}}\frac{a}{h_u}=e^{y+\frac12a^2y^2}=
\sum_{n=0}^{\infty}\frac{y^n}{n!}\sum\Sb\pi\in Inv(S_n)\endSb a^{2c_2(\pi)}.$$ \rm
\bf Example 3.7 \rm Put $q=t=\frac1{z}, x=yz$  in (1), so that
$(1-yz)^{-\frac1z}=\sum_n\frac{y^n}{n!}\sum_{\pi\in S_n}z^{n-\kappa(\pi)}$.
Now, take the limit $z\rightarrow 0$ and note $n=\kappa(\pi)$ iff $\pi=1$. The outcome relates to the classical result
$$\sum_n\frac{y^n}{n!^2}\sum_{\lambda\vdash n}f_{\lambda}^2=e^y=\sum_n\frac{y^n}{n!}\sum\Sb\pi\in S_n\\ \pi=1\endSb 1.$$ \rm
\bf Example 3.8 \rm Put $q=t^2=\frac1{z^2}, x=yz$ in (1), then take $z\rightarrow 0$. Observe again that
$n=\Cal{O}(\pi)+2\Cal{E}(\pi)$ iff $\pi^2=1$. The result relates to another classical formula
$$\sum_n\frac{y^n}{n!}\sum_{\lambda\vdash n}f_{\lambda}=e^{y+\frac12y^2}=
\sum_n\frac{y^n}{n!}\sum_{\pi\in Inv(S_n)} 1.$$ \rm

\smallskip
\noindent
\bf 4 OF A DIFFERENT BREED \rm
\bigskip
\noindent
Often, we see $tr(\pi)=c_1(\pi)$ as an \bf exponent. \rm What if it serves as a \bf base \rm instead?
\bigskip
\noindent
\bf Conjecture 4.1 \it We have the generating function
$$\sum_{n=0}^{\infty}\frac{z^n}{n!}\sum\Sb\pi\in Inv(S_n)\endSb tr(\pi)^k=\widehat{B}_k(z)e^{z+\frac12z^2}.$$ \rm
\bf Conjecture 4.2 \it We have the generating function
$$\sum_{n=0}^{\infty}\frac{x^n}{n!}\sum\Sb\pi\in S_n\\ \Cal{E}(\pi)=0\endSb 2^{\kappa(\pi)}\widehat{B}_{\kappa(\pi)}(z)
=e^{\frac{2zx}{1-x}}.$$
\bf Conjecture 4.3 \it Interesting consequences are
$$\align
\frac{\sum_{\pi\in Inv(S_n)}tr(\pi)^k}{n!}&=\sum_{j=0}^kS(k,j)\frac{\#Inv(S_{n-j})}{(n-j)!}, \\
\sum\Sb\pi\in S_n\\ \Cal{E}(\pi)=0\endSb 2^{\kappa(\pi)}B_{\kappa(\pi)}(z)&=
\sum\Sb\pi\in S_n\endSb (2z)^{\kappa(\pi)}\prod_{j=1}^nj^{c_j(\pi)}.\endalign$$ \rm
\bf Remark. \rm The specialization at $z=1$ implies
$$\sum\Sb\lambda\vdash n\\ \Cal{E}(\lambda)=0\endSb \frac{2^{L(\lambda)}B_{L(\lambda)}}{1^{m_1}3^{m_3}5^{m_5}\cdots m_1!m_3!m_5!\cdots}=\sum_{\lambda\vdash n}\frac{2^{L(\lambda)}}{m_1!m_2!m_3!\cdots}.$$

\bigskip
\noindent
\bf 5 CHIRALS in SUPERSYMMETRY GAUGE THEORY \rm
\bigskip
\noindent
We start with a lemma.
\smallskip
\noindent
\bf Lemma 5.1 \it Using the \bf bra-ket \it notation, the partition function $Z_{U(1)}$ factorizes as
$$\sum_{k=0}^r\alpha_k^r\langle \text{Tr $\psi^{2k}$}\rangle Z_{U(1)}=
\left[2\binom{2r-1}{r-1}^2q^r+\frac{r^2}2\binom{2r}r\binom{2r-4}{r-2}q^{r-1}\right]e^q.$$ \rm
\bf Proof. \rm In $U(1)$ supersymmertic gauge theory, there are chiral observables (Hermitian operators) $\text{Tr $\psi^k$}$, where $\psi$ is the Higgs scalar field in the adjoint representation. The scalar correlators compute as $1$-point functions $\langle \text{Tr $\psi^{2k}$}\rangle$, in the parameter of instanton expansion $q$:
$$\langle \text{Tr $\psi^{2k}$}\rangle=\frac1{Z_{U(1)}}\sum_{n=0}^{\infty}\sum_{\lambda\vdash n}
\frac{ \text{Tr $\psi_{\lambda}^{2k}$}}{H_{\lambda}^2}q^n,$$
where $H_{\lambda}^{-2}=\prod_{u\in\lambda}h_u^{-2}$ is the Plancherel measure $\mu(\lambda)^2$ on the space of (random) partitions (it can be regarded as the discretization of the Vandermonde measure on random matrices);
$\text{Tr $\psi_{\lambda}^{2k}$}=\sum_{j=0}^{k-1}2\binom{2k}{2j}\sum_{u\in\lambda}h_u^{2j}$; and $Z_{U(1)}$ is the partition function of $U(1)$ gauge theory
$$Z_{U(1)}=\sum_{n=0}^{\infty}\sum_{\lambda\vdash n}\frac1{H_{\lambda}^2}q^n.$$
If we denote $\Omega_j(n)=\sum_{\lambda\vdash n}\mu(\lambda)^2\sum_{u\in\lambda}h_u^{2j}$, then
$$\langle \text{Tr $\psi^{2k}$}\rangle Z_{U(1)}=\sum_{n\geq0}\sum_{j=0}^{k-1}2\binom{2k}{2j}\Omega_j(n) q^n.$$
Denote $\Cal{P}_{2m}(x):=\prod_{j=0}^{m-1}(x^2-j^2)=\sum_{j=1}^m\alpha_j^mx^{2j}$ with
$\alpha_j^m=(-1)^{m-j}e_{m-j}\left(1^2,2^2,\dots,(m-1)^2\right)$.
Summing over all instantons, results in
$$\align
\sum_{k=0}^r\alpha_k^r\langle \text{Tr $\psi^{2k}$}\rangle Z_{U(1)}&=
\sum_{n\geq0}q^n\sum_{\lambda\vdash n}\sum_{k=0}^r\alpha_k^r
\frac{\sum_{u\in\lambda}(h_u+1)^{2k}-2h_u^{2k}+(h_u-1)^{2k}}{\prod_{u\in\lambda}h_u^2}\\
&=\sum_{n\geq0}q^n\sum_{\lambda\vdash n}
\frac{\sum_{u\in\lambda}\Cal{P}_{2r}(h_u+1)-2\Cal{P}_{2r}(h_u)+\Cal{P}_{2r}(h_u-1)}
{\prod_{u\in\lambda}h_u^2}\\
&=2r(2r-1)
\sum_{n\geq0}q^n\sum_{\lambda\vdash n}\sum_{u\in\lambda}\frac{\prod_{j=0}^{r-2}(h_u^2-j^2)}{H_{\lambda}^2},
\endalign$$
where we have utilized $\Cal{P}_{2r}(x)=(x)^{r}\cdot(x)_{r}$ and the discrete "derivative" relation
$$\Delta^2\Cal{P}_{2r}(x):=\Cal{P}_{2r}(x+1)-2\Cal{P}_{2r}(x)+\Cal{P}_{2r}(x-1)=2r(2r-1)\Cal{P}_{2r-2}(x).\qquad
\square$$
\bf Corollary 5.2 \it We have the combinatorial identity
$$\frac1{n!}\sum_{\lambda\vdash n}f_{\lambda}^2\sum_{u\in\lambda}\prod_{\pmb{j=0}}^{r-1}(h_u^2-j^2)=
\binom{2r+1}{r+1}^2\binom{n}{r+1}\frac{r!}{2r+1}+
\binom{2r+2}{r+1}\binom{2r-2}{r-1}\binom{n}r\frac{(r+1)!}{8r+4}.$$ \rm
\bf Problem 5.3 \rm Find a combinatorial proof for Cor. 5.2.

\bigskip
\noindent
\bf 6 SYMPLECTIC and ORTHOGONAL SEMISTANDARD YOUNG DIAGRAMS \rm
\bigskip
\noindent
Recall the hook-length $h(i,j)=\lambda_i+\lambda_j^t-i-j+1$ and content $c(i,j)=j-i$ of a cell $u=(i,j)$ of a Young diagram of shape $\lambda$.
Define the \it symplectic content \rm for symplectic
diagram as
$$c_{sp}(i,j)=\cases \lambda_i+\lambda_j-i-j+2 \qquad \text{if $i>j$} \\
i+j-\lambda_i^t-\lambda_j^t \qquad \qquad \text{if $i\leq j$},\endcases$$
and the \it orthogonal content \rm for orthogonal Young diagram as
$$c_{O}(i,j)=\cases \lambda_i+\lambda_j-i-j \qquad \qquad \text{if $i\geq j$} \\
i+j-\lambda_i^t-\lambda_j^t-2 \qquad \text{if $i<j$}.\endcases$$
The \it symplectic Schur function \rm corresponding to $\lambda$ is computed by the determinantal formula [5]
$$sp_{\lambda}(x_1,\dots,x_n)=\frac{\det\left(x_j^{\lambda_i+n-i+1}-x_j^{-\lambda_i-n+i-1}\right)}
{\det\left(x_j^{n-i+1}-x_j^{-n+i-1}\right)}.$$
Interesting results such as Stanley's [11, Cor. 7.21.4] hook-content formula
$s_{\lambda}(1^n)=\prod_{u\in\lambda}\frac{n+c_u}{h_u}$ find their symplectic counterpart
[2, Cor. 9]
$$sp_{\lambda}(1^n)=\prod_{u\in\lambda}\frac{2n+c_{sp}(u)}{h(u)}.$$
\bf Proposition 6.1 \it We have the identity
$$\sum_{\lambda\vdash n}\prod_{u\in\lambda}\frac{c_{sp}(u)}{h(u)}=
\sum_{\lambda\vdash n}\prod_{u\in\lambda}\frac{c_{O}(u)}{h(u)}.$$ \rm
\bf Proof. \rm This is immediate from the definitions of $c_{sp}, c_O$ and the fact that
$\sum_{\lambda\vdash n}=\sum_{\lambda^t\vdash n}$. $\square$

$$ $$

\pagebreak

\bigskip
\noindent
\bf Conjecture 6.2(a) \it We have the generating function
$$\sum_{n\geq0}x^n\sum_{\lambda\vdash n}\prod_{u\in\lambda}\frac{t+c_{sp}(u)}{h(u)}=
\prod_{j\geq1}\frac{(1-x^{8j})^{\binom{t+1}2}}{(1-x^{8j-2})^{\binom{t+1}2-1}}
\left(\frac{1-x^{4j-1}}{1-x^{4j-3}}\right)^t
\left(\frac{1-x^{8j-4}}{1-x^{8j-6}}\right)^{\binom{t}2-1}.$$ \rm
\bf Conjecture 6.2(b) \it We have the generating function
$$\sum_{n\geq0}x^n\sum_{\lambda\vdash n}\prod_{u\in\lambda}\frac{t+c_{O}(u)}{h(u)}=
\prod_{j\geq1}\frac{(1-x^{8j})^{\binom{t}2}}{(1-x^{8j-6})^{\binom{t}2-1}}
\left(\frac{1-x^{4j-1}}{1-x^{4j-3}}\right)^t
\left(\frac{1-x^{8j-4}}{1-x^{8j-2}}\right)^{\binom{t+1}2-1}.$$ \rm
\bf Conjecture 6.2(c) \it The specialization at $t=0$ gives
$$\sum_{n\geq0}x^n\sum_{\lambda\vdash n}\prod_{u\in\lambda}\frac{c_{sp}(u)}{h(u)}=
\sum_{n\geq0}x^n\sum_{\lambda\vdash n}\prod_{u\in\lambda}\frac{c_{O}(u)}{h(u)}=
\prod_{j\geq1}\frac1{1+x^{4j-2}}.$$ \rm
\bf Conjecture 6.3(a) \it We have a symplectic/orthogonal "content Nekrasov-Okounkov formula":
$$\sum_{n\geq0}x^n\sum_{\lambda\vdash n}\prod_{u\in\lambda}\frac{t+c_{sp}^2(u)}{h^2(u)}=
\sum_{n\geq0}x^n\sum_{\lambda\vdash n}\prod_{u\in\lambda}\frac{t+c_{O}^2(u)}{h^2(u)}=
\prod\Sb j\geq1\endSb\frac1{1-x^{4j-2}}\prod_{j\geq1}\frac1{(1-x^j)^t}.$$ \rm
\bf Conjecture 6.3(b) \it The specialization at $t=0$ gives
$$\sum_{n\geq0}x^n\sum_{\lambda\vdash n}\prod_{u\in\lambda}\frac{c_{sp}^2(u)}{h^2(u)}=
\sum_{n\geq0}x^n\sum_{\lambda\vdash n}\prod_{u\in\lambda}\frac{c_{O}^2(u)}{h^2(u)}=
\prod_{j\geq1}\frac1{1-x^{4j-2}}.$$ \rm
\bf Conjecture 6.3(c) \it Let $\delta_n=(-1)^{\binom{n}2}$. Combining conjectures 6.2(c) and 6.3(b),
we have the identity
$$\delta_n\sum_{\lambda\vdash n}f_{\lambda}^2\prod_{u\in\lambda}c_{sp}^2(u)=n!
\sum_{\lambda\vdash n}f_{\lambda}\prod_{u\in\lambda}c_{sp}(u).$$ \rm
\bf Problem 6.4 \it In the classical Schur and elementary symmetric functions, apply the "flattening conditions" of replacing all powers $y_i^m\rightarrow y_i$. Then, find a \bf direct \rm proof for
$$\sum_{\lambda\vdash k}f_{\lambda}s_{\lambda}(y_1,\dots,y_n)=
\sum_{j} j! S(k,j) e_j(y_1,\dots,y_n).$$ \rm
\bigskip
\noindent
\bf 7 ENUMERATION under GROUP ACTION \rm
\bigskip
\noindent
P\'olya contributed one of his pioneering works to this field; for a brief recall see [10, Sec. 7.24], and references therein. At present, we confine the discussion to the \it cycle index polynomial \rm [11, Example 5.2.10] of the symmetric group $S_n$:
$$Z(S_n)=\frac1{n!} \sum_{\pi\in S_n}Z(\pi,\pmb{t});$$
where $Z(\pi,\pmb{t})=t_1^{c_1(\pi)}t_2^{c_2(\pi)}\cdots t_n^{c_n(\pi)}$. With this on the background, we opt to propose the following determinantal expansion.
\bigskip
\noindent
\bf Problem 7.1 \it Given a smooth differential operator $L$, linearize it by local charts so that L acts as a
matrix, say $n$-by-$n$. Denote $t_i=tr(L^i)$ as the trace of a matrix. Then, we have
$$\det(L)=\frac1{n!} \sum_{\pi\in S_n}(-1)^{\kappa(\pi)-1}Z(\pi,\pmb{t}).\tag2$$
It's probably not hard to employ symmetric functions in proving (2). However, the above determinant behaves errily similar to a cycle index formula. So, find an interpretation of (2) as some appropriate group action. \rm

\pagebreak
\bigskip
\noindent
\bf 8 QUANTUM VERSIONS of HAN'S IDENTITIES \rm
\bigskip
\noindent
Given a complex number $\alpha$, G. Han [7] has shown that
$$\sum_{n\geq1}x^n\sum_{\lambda\vdash n}\sum_{u\in\lambda}h_u^{\alpha}=
\prod_{m\geq1}\frac1{1-x^m}\times\sum_{k\geq1}\frac{x^kk^{\alpha+1}}{1-x^k}=
\sum_{n\geq1}x^n\sum_{\lambda\vdash n}\sum_{i}\lambda_i^{\alpha+1}.\tag3$$
Now, we state the following generalizations (compute $\frac{d}{dq}$ and set $q=1$ to recover (3) above).
\smallskip
\noindent
\bf Conjecture 8.1 \it For any complex number $\alpha$ and a free parameter $q$, we have
$$\align
\sum_{n\geq1}x^n\sum_{\lambda\vdash n}\sum_{u\in\lambda}q^{h_u^{\alpha}}&=
\prod_{m\geq1}\frac1{1-x^m}\times\sum_{k\geq1}\frac{kx^kq^{k^{\alpha}}}{1-x^k}\\
\sum_{n\geq1}x^n\sum_{\lambda\vdash n}\sum_{i}q^{\lambda_i^{\alpha}}&=
\prod_{m\geq1}\frac1{1-x^m}\times\sum_{k\geq1}\frac{x^kq^{k^{\alpha}}}{1-x^k}.\endalign$$ \rm
\bf Proposition 8.2 \it Let $\zeta(s)$ be the Riemann zeta function. We have the following variants:
$$\align
\sum_{n\geq1}x^n\sum_{\lambda\vdash n}\sum_{i}\lambda_i &=
\prod_{m\geq1}\frac1{1-x^m}\times\sum_{k\geq1}\frac{x^k}{(1-x^k)^2} \\
\sum_{n\geq1}x^n\sum_{\lambda\vdash n}\sum_{i}q^{\lambda_i} &=
\prod_{m\geq1}\frac1{1-x^m}\times\sum_{k\geq1}\frac{qx^k}{1-qx^k}\\
\sum_{\lambda\vdash n}\sum\Sb i\\  \lambda_i\neq1\endSb\frac1{\lambda_i-1}&=
\sum_{\lambda\vdash n}\sum\Sb u\in\lambda\\ h_u\neq1\endSb\frac1{h_u(h_u-1)}\\
\sum_{\lambda\vdash n}\sum\Sb i\\  \lambda_i\neq1\endSb\lambda_i\zeta(\lambda_i)&=
\sum_{\lambda\vdash n}\sum\Sb u\in\lambda\\ h_u\neq1\endSb\zeta(h_u).\endalign$$ \rm
\bf Proof. \rm This is straightforward: $i\cdot \#\{j\vert \lambda_j=i;\lambda\vdash n\}=
\#\{u\vert h_u=i; \lambda\vdash n\}$. $\square$
\bigskip
\noindent
\bf 9 COUNTING SQUARES IN YOUNG DIAGRAMS \rm
\bigskip
\noindent
This section has benefited from contributions by R. Stanley.
First, we state the following celebrated result which is somehow linked to counting squares
in the Young diagram of a partition.
\smallskip
\noindent
\bf Theorem (First Roger-Ramanjuan identities) \it Given $n\in\Bbb{N}$, define the two sets
$$\align
A_n&=\{\lambda\vdash n: \lambda_i-\lambda_{i+1}\geq2, \forall i\}=\text{partitions with parts differing at least by $2$}, \\
B_n&=\{\lambda\vdash n: \lambda_i\equiv 1,4\mod5, \forall i\}=\text{partitions with parts congruent to $1$ or $4$ mod 5}. \endalign$$
Then, $\#A_n=\#B_n$. Equivalently,
$$1+\sum_{k\geq1}\frac{x^{k^2}}{(x;x)_k}=\prod_{j\geq1}\frac1{(1-x^{5j-1})(1-x^{5j-4})}.$$ \rm
\bf Proposition 9.1 \it For the set $A_n$, we have
$$\sum_{n\geq1}x^n\sum_{\mu\in A_n}q^{h(1,1)}=
\sum_{k\geq1}\frac{x^{k^2}q^{3k-2}}{(qx;x)_k}.$$ \rm
\bf Proof. \rm If $\mu\in A_n$, then $n=\mu_1+\cdots+\mu_k=(y_1+2k-1)+\cdots+(y_{k-1}+3)+(y_k+1)$ for some
$y\vdash n-k^2$. Notice $h(1,1)=\mu_1+k-1=y_1+3k-2$. So, the generating function expands as:
$$\align \sum_{n\geq1}x^n\sum_{\lambda\in A_n}q^{h(1,1)}
&=\sum_{k\geq1}x^{k^2}q^{3k-2}\sum_n\sum\Sb y\vdash n-k^2\endSb x^{y_1+\cdots+y_k}q^{y_1} \\
&=\sum_{k\geq1}x^{k^2}q^{3k-2}\sum_{0\leq z_k\leq\cdots\leq z_1}x^{z_1+\cdots+z_k}q^{z_1} \\
&=\sum_{k\geq1}\frac{x^{k^2}q^{3k-2}}{1-xq}\sum\Sb 0\leq z_k\leq\cdots\leq z_2\endSb x^{2z_2}x^{z_3+\cdots+z_k}q^{2z_2}\\
&=\dots\\
&=\sum_{k\geq1}\frac{x^{k^2}q^{3k-2}}{(1-xq)\cdots(1-x^kq)}.\qquad\square \endalign$$
\bf Proposition 9.2 \it Let $[q^m]F(q)=$ coefficient of $q^m$ in the Taylor series for $F(q)$. Then,
$$\sum_{n\geq0}x^n\sum_{\lambda\vdash n}q^{h(1,1)}=
\sum_{k\geq0}\frac{x^{k(k+1)}q^{2k}}{(qx;x)_k(qx;x)_{k+1}}
=\sum_{n\geq0}x^n\sum_{i,j}q^j\left[q^{n-j}\right]\pmb{\binom{j-1}i}.$$ \rm
\bf Notations. \rm Write the infinite vector $\pmb{N}:=\langle1,2,3,\dots\rangle$ and "$\bullet$" for the usual dot product.
Given a partition $\lambda\vdash n$, pad up infinitely many zeros at the end to denote $\lambda=\langle
\lambda_1,\lambda_2,\dots\rangle$ and its diagonal-hooks by $\pmb{h}_{\lambda}=\langle h(1,1),h(2,2),\dots\rangle$ to be treated as vectors.

\bigskip
\noindent
\bf Lemma 9.3 \it Denote the set of partitions of $n$ by $\Cal{P}_n:=\{\lambda: \lambda\vdash n\}$ and define the \bf multiset \it
of hooks as $C_n:=\{\pmb{h}_{\lambda}:\lambda\vdash n\}$. Then,

(a) the map $\Psi_n:\Cal{P}_n\rightarrow C_n$ given by
$\Psi_n(\lambda)=\pmb{h}_{\lambda}$ is onto $A_n$.

(b) In particular, $\pmb{h}_{\lambda}\vdash n$; i.e. each partition $\lambda\vdash n$
generates another partition $\mu=\pmb{h}_{\lambda}$ whose

parts differ by at least $2$. \rm
\bigskip
\noindent
\bf Proof. \rm First, observe that if $\lambda\vdash n$ then $\sum_{i\geq1}h(i,i)=n$ and $h(i,i)-h(i+1,i+1)\geq2$.
Thus $\pmb{h}_{\lambda}\in A_n$. Of course, by definition $\pmb{h}_{\lambda}\in C_n$. Although $A_n\subset C_n$, when the latter is
viewed as a multiset, once repetitions are removed from $C_n$ then both $A_n$ and $C_n$ become exactly equal as \bf sets. \rm
Therefore, as sets (not multisets) $\#A_n=\#B_n=\#C_n$. The proof is complete. $\square$
\bigskip
\noindent
\bf Definition 9.4 \rm Given a partition $\lambda$, define
\smallskip
(i) $a(\lambda)$ to be the \it number of squares \rm (of all sizes) in the Young diagram to $\lambda$;

(ii) the sequence $f(n)=\sum_{\lambda\vdash n}a(\lambda)$;

(iii)  the polynomials $F_n(q)=\sum_{\lambda\vdash n}q^{a(\lambda)}$.

$$ $$

\pagebreak

\bigskip
\noindent
\bf Theorem 9.5 \it Let $i\wedge j=\min\{i,j\}$. We have the following properties:

\smallskip
(i) $a(\lambda)=\pmb{N}\bullet\pmb{h}_{\lambda}=\sum_{(i,j)\in\lambda}i\wedge j$.

(ii) $F_n(1)=p(n)$ the number of partitions of $n$,  and $\frac{d}{dq}F_n(1)=f(n)$.

(iii) the generating function
$$\sum_{n\geq0}F_n(q)x^n=
\sum_{k\geq0}\frac{x^{k^2}q^{\frac{k(k+1)(2k+1)}6}}{\prod_{j=1}^k\left(1-x^jq^{\frac{j(j+1)}2}\right)^2}.$$ \rm
\bf Proof. \rm (i) Identify squares by their bottom-right corners (BLC) and mark each cell $(i,j)\in\lambda$ by the number of
squares for which it is a $BLC$. This easily enumerates as $\#BLC(i,j)=i\wedge j$, hence both claims $\sum_{(i,j)}i\wedge j$ and $\pmb{N}\bullet\pmb{h}_{\lambda}$ are now immediate. (ii) Obvious. (iii) Begin with the well-known identity $\sum_n\sum_{\lambda\vdash n}x^n=\prod_j\frac1{1-x^j}=\sum_k\frac{x^{k^2}}{\prod_{j=1}^k(1-x^j)^2}$. On the basis of Lemma 9.3, expand the left-hand side in the form (denote $\Cal{P}=\bigcup_n\Cal{P}_n=$ set of all partitions, $\sigma_k=(2k-1,\dots,2,1)$, use $\mu_j=y_j+2k-2j+1$ to follow the proof of Proposition 9.1):
$$\align\sum_n\sum_{\lambda\vdash n}x^n=
\sum_n\sum_{\mu\in A_n}\sum\Sb\lambda\vdash n\\ \pmb{h}_{\lambda}=\mu\endSb x^{\mu_1+\mu_2+\cdots}&=
\sum_kx^{k^2}\sum_n\sum\Sb 0\leq y_k\leq\cdots\leq y_1\endSb\sum\Sb \lambda\vdash n\\ \pmb{h}_{\lambda}=y+\sigma_k\endSb x^{y_1+\cdots+y_k} \\
&=\sum_kx^{k^2}\sum\Sb 0\leq y_k\leq\cdots\leq y_1\endSb\sum\Sb \lambda\in\Cal{P}\\ \pmb{h}_{\lambda}=y+\sigma_k\endSb x^{y_1+\cdots+y_k}.
\tag*\endalign$$
In light of this, we obtain
$$\align\sum_{n\geq0}\sum_{\lambda\vdash n}x^nq^{a(\lambda)}
=&\sum_{n\geq0}\sum_{\mu\in A_n}\sum\Sb\lambda\vdash n\\ \pmb{h}_{\lambda}=\mu\endSb x^{\mu_1+\mu_2+\cdots}q^{\mu_1+2\mu+3\mu_3+\cdots} \\
=&\sum_kx^{k^2}q^{\frac{k(k+1)(2k+1)}6}
\sum\Sb 0\leq y_k\leq\cdots\leq y_1\endSb\sum\Sb \lambda\in\Cal{P}\\ \pmb{h}_{\lambda}=y+\sigma_k\endSb
 (xq)^{y_1}(xq^2)^{y_2}\cdots (xq^k)^{y_k}. \tag**
\endalign$$
Since (*) equates to $\sum_k\frac{x^{k^2}}{\prod_{j=1}^k(1-x^j)^2}$, it is rather evident how (**)
leads to the assertion of part (iii). The proof is complete. $\square$

\bigskip
\noindent
\bf Definition  9.6 \rm We say two partitions $\lambda, \overline{\lambda}\vdash n$ are related,
$\lambda\sim\overline{\lambda}$, provided that $\pmb{N}\bullet\pmb{h}_\lambda=\pmb{N}\bullet\pmb{h}_{\overline{\lambda}}$.
Recall that $\pmb{h}_{\lambda}$ and $\pmb{h}_{\overline{\lambda}}$ belong to $A_n$. Evidently, too, $\sim$ is an
equivalence relation on $\Cal{P}_n$, where the set $A_n$ plays an internal role.
\bigskip
\noindent
\bf Corollary 9.7 \it For each $j\in\Bbb{N}$, denote the equivalent classes
$D_j(n)=\{\lambda\in\Cal{P}_n: \pmb{N}\bullet\pmb{h}_{\lambda}=j\}$. Then,

\smallskip
(i) $\bigcup_{j=n}^{\infty} D_j(n)$ is a (finite) disjoint partition of the set of partitions $\Cal{P}_n$.

(ii) $\#D_j(n)=$ the coefficient of $q^j$ in the polynomial $F_n(q)$.

(iii) write each $\mu\vdash j$ as $\mu=1^{m_1}2^{m_2}\cdots$, then the number of partitions of $j$ such that

$j=1^{h(1,1)}2^{h(2,2)}\cdots$ where the $h$'s are hooks for $\lambda\vdash n$ equals $[q^j]F_n(q)$. \rm
\bigskip
\noindent
\bf Question 9.8 \it By definition, $f(n)$ enumerates the total number of squares in all partitions of $n$. What
else can be said about $f(n)$? For instance, is there some representation-theoretic meaning associated to it? \rm

\pagebreak

\bigskip
\noindent
\bf 10 A THEOREM OF FROBENIUS \rm
\bigskip
\noindent
It's fair to say, F. G. Frobenius created representation theory of finite groups after receiving a letter from R. Dedekind.
This letter contained the following procedure: let $G$ be a finite group and $\det(G)(X)$ be the \it group determinant \rm of $G$ (formed by its multiplication table, regarding the group elements as \it commuting \rm
indeterminates). Then, a theorem of Frobenius factors $\det(G)(X)$ into a
product of irreducible polynomials $P_j(X)$ (over the complex numbers) with
multiplicity equaling the degree of $P_j(X)$; that is,
$$\det(G)(X)=\prod_{j=1}^rP_j(X)^{\deg P_j}$$
where $r$ is the number of conjugacy classes of $G$.
\bigskip
\noindent
\bf Question 10.1. \rm Is the converse true? Given just the multiplication table for a
Latin square, compute the determinant, and assume that it factors into
irreducible polynomials with multiplicity equal to the degree of the corresponding polynomial as in the theorem of Frobenius. We know that if it does not factor as in the theorem then the
Latin square is not that of a group. If it does factor well, can we
conclude that the table, up to a permutation of its rows, is the
multiplication table of a group?

\bigskip
\noindent
\bf 11 ON $t$-CORE PARTITIONS \rm
\bigskip
\noindent
For a positive integer $s$, we say a partition $\lambda$ is a \it $t$-core \rm if $\lambda$ has no hook of length equal to $t$. Similarly, $\lambda$ is an \it $(r,s,t)$-core \rm if it simultaneously is an $r$- and an $s$- and a $t$-core. For related results on $(s,t)$-cores, see $[1]$ and [12] and references therein.
\bigskip
\noindent
\bf Conjecture 11.1 \it Let $C_k=\frac1{k+1}\binom{2k}k$ be the Catalan number. The number $f(s)$ of $(s,s+1,s+2)$-cores equals to
$$f(s)=\sum_{k\geq0}\binom{s}{2k}C_k.$$ 
Generating function: $\sum f(s)x^s=\frac{1\vert}{\vert 1-x}-\frac{x^2\vert}{\vert 1-x}-
\frac{x^2\vert}{\vert 1-x}-\cdots=\frac{1-x-\sqrt{1-2x-3x^2}}{2x^2}$. \rm
\smallskip
\noindent
\bf Conjecture 11.2 \it The size $g(s)$ of largest $(s,s+1,s+2)$-core partition (based on the parity of $s$) equals to
$$g(s)=\cases \qquad m\binom{m+1}3, \qquad \qquad \qquad \text{if $s=2m-1$} \\
(m+1)\binom{m+1}3+\binom{m+2}3, \qquad \text{if $s=2m$}. \endcases$$ \rm
\bf Conjecture 11.3 \it The total sum $h(s)$ of the sizes of all $(s,s+1,s+2)$-cores equals to
$$h(s)=\sum_{j=0}^{s-2}\binom{j+3}3\sum_{i=0}^{\lfloor{j/2}\rfloor}\binom{j}{2i}C_i.$$
Equivalently, the average size of an $(s,s+1,s+2)$-core is $\frac{h(s)}{f(s)}$. \rm

\pagebreak

\smallskip
\noindent
The results are even more general as stated next.
\smallskip
\noindent
\bf Conjecture 11.4 \it Fix $p\geq2$, Then, the number $F_p(s)$ of $(s,s+1,\dots,s+p)$-core partitions is computed by the generalized Catalan numbers, and is given recursively by $F_p(0)=1$ so that
$$F_p(s+1)=F_p(s)+\sum_{k=p-2}^{s-1}F_p(k)F_p(s-1-k).$$
Convention: empty sums are zero. For $p\geq1$, there is also the generating function
$$\sum_{s\geq0}F_p(s)x^s=\frac{A_p(x)-\sqrt{A_p(x)^2-4x^2}}{2x^2}, \qquad\text{where} \qquad
A_p(x)=1-x+\frac{x^2-x^p}{1-x}.$$
Equivalently, we have the continued fraction expansion
$$\align \sum_{s\geq0}F_p(s)x^s
&=\frac{1\vert}{\vert A_p(x)}-\frac{x^2\vert}{\vert A_p(x)}-\frac{x^2\vert}{\vert A_p(x)}-\cdots. \endalign$$ \rm

\smallskip
\noindent
Encouraged by this success, we offer the following generalized conjecture (joint with Emily Leven) regarding partitions whose cores line up in arithmetic progression.

\smallskip
\noindent
\bf Conjecture 11.5. \it Let $s$ and $d$ be two relatively prime positive integers. Then the number of $(s,s+d,s+2d)$-core partitions is given by
$$\sum_{k=0}^{\lfloor\frac{s}2\rfloor}\binom{s+d-1}{2k+d-1}\binom{2k+d}k\frac1{2k+d}
=\frac1{s+d}\sum_{k=0}^{\lfloor\frac{s}2\rfloor}\binom{s+d}{k,k+d,s-2k}.$$ \rm

\noindent
The special case $d=1$ simply becomes the \it Motzkin numbers. \rm This connection perpetuates to the case $d=2$ as stated below. 

\smallskip
\noindent
\bf Problem 11.6. \rm The \it Motzkin triangle \rm $T(n,k)$ is built according to the rules:

(1) $T(n,0)=1$;

(2) $T(n,k)=0$ if $k<0$ or $k>n$;

(3) $T(n,k)=T(n-1,k-2)+T(n-1,k-1)+T(n-1,k)$.

\smallskip
\noindent
Prove the identity
$$\sum_{k=0}^nT(n,k)T(n,k+1)=\sum_{k=0}^n\binom{2n}{2k+1}\binom{2k+1}k\frac1{k+2}.$$ 
\bf Conjecture 11.7. \it The number of $(s,s+1)$-core partitions $\lambda$ with parts that are multiples of $\ell$ equals
$$\frac{(s\mod \ell)+1}{s+1}\binom{(\ell+1)\lfloor\frac{s}{\ell}\rfloor+(s\mod\ell)}{\lfloor\frac{s}{\ell}\rfloor}
=\frac{s+1-\ell\lfloor\frac{s}{\ell}\rfloor}{s+1}\binom{s+\lfloor\frac{s}{\ell}\rfloor}s.$$
In particular, if $\ell=1$ then this count is the Catalan number $C_s=\frac1{s+1}\binom{2s}s$. \rm
\bigskip
\noindent
\bf Easy exercise 11.8. \it The generating functions for the number $c_2(n)$ of $2$-core and $c_3(n)$ of $3$-core partitions $\lambda\vdash n$ with distinct parts are, respectively,
$$\align \sum_{n\geq0}c_2(n)q^n&=\prod_{k\geq1}(1-q^{2k})(1+q^k)=\sum_{n\geq0}q^{\binom{n+1}2} \qquad \text{and} \\
\sum_{n\geq0}c_3(n)q^n&=\frac12\prod_{k\geq1}(1-q^{2k})(1+q^{2k-1})^2+
\prod_{k\geq1}(1-q^{2k})(1+q^{2k})^2-\frac12=\sum_{n\geq0}q^{n^2}+\sum_{n\geq1}q^{2\binom{n+1}2}. \endalign$$ \rm
\bf Conjecture 11.9. \it For $(s,s+1)$-core partitions with distinct parts, we have
\smallskip
(a) the number of such partitions equals the Fibonacci number $F_{s+1}$;
\smallskip
(b) the largest size of such partition is $\left\lfloor\frac13\binom{s+1}2\right\rfloor$;
\smallskip
(c) there are $\frac{3-(-1)^{s\,\,\text{mod}\,\, 3}}2$ such partitions of maximal size;
\smallskip
(d) the total number of these partitions and the average sizes are, respectively, given by
$$\sum_{i+j+k=s+1}F_iF_jF_k \qquad \text{and} \qquad \sum_{i+j+k=s+1}\frac{F_iF_jF_k}{F_{s+1}}.$$ \rm
\bf Conjecture 11.10. \it For $s$ odd, the number of $(s,s+2)$-core partitions with distinct parts is $2^{s-1}$. \rm
\bigskip
\noindent
\bf Easy exercise 11.11. \it The generating function for the number $\text{oc}_3(n)$ of $3$-core partitions $\lambda\vdash n$ with odd parts is given by
$$\sum_{n\geq1}\text{oc}_3(n)q^n=\frac{1+q}2\left[\prod_{k\geq1}(1-q^{2k})(1+q^{2k-1})^2-1\right]
=(1+q)\sum_{n\geq1}q^{n^2}.$$ \rm

\bigskip
\noindent
\bf 12 ON COLORED OVERPARTITIONS \rm
\bigskip
\noindent
In the literature an \it overpartiton \rm $\lambda\vdash n$ is defined as an ordinary integer partition where the last occurrence of each distinct part in $\lambda$ may be overlined. The number of overpartitions of $n$ is customerly denoted by $\overline{p}(n)$. For example, the $8$ overpartitions of $3$ are:
$$3,\, \overline{3}, \, 2+1, \, \overline{2}+1, \, 2+\overline{1}, \, \overline{2}+\overline{1}, \, 1+1+1, \, \overline{1}+1+1.$$
Next, we introduce a natural generalization of this concept.
\bigskip
\noindent
\bf Definition. \rm A \it colored overpartition \rm (COP) is a partition where the last occurrence of each distinct number may receive any one of $c$ colors. The number of such partitions of $n$ we denote by $\overline{p}_c(n)$.
\bigskip
\noindent
\bf Note. \rm The number of: (i) ordinary partitions $p(n)=\overline{p}_1(n)$; (ii)  overpartitions $\overline{p}(n)=\overline{p}_2(n)$. 
\bigskip
\noindent
It is not a big effort to arrive at the following which incidentally extend some results of [13].
\smallskip
\noindent
\bf Theorem 12.1 \it A generating function for the number of COP fitting inside an $M\times N$ rectangle is
$$\sum_{k=0}^{\min\{M,N\}}(c-1)^kq^{\binom{k+1}2}\frac{(q)_{M+N-k}}{(q)_k(q)_{M-k}(q)_{N-k}}.$$
The limiting case (both $M, N\rightarrow\infty$) provides a generating function for the number of COP
$$\sum_{n\geq0}\overline{p}_c(n)q^n=\frac1{(q;q)_{\infty}}\sum_{k=0}^{\infty}
\frac{(c-1)^kq^{\binom{k+1}2}}{(q;q)_k}=\frac{((1-c)q;q)_{\infty}}{(q;q)_{\infty}}.$$ \rm
We are set to propose an overarching congruence valid for all colored overpartitions $\overline{p}_c(n)$, collectively.
\smallskip
\noindent
\bf Conjecture 12.2 \it Let $c\in\Bbb{N}$. We have the following congruence relations
$$\overline{p}_c\left(2^cn+2^{c-1}\right)\equiv 0\mod c^2,\qquad \text{for all $n\geq0$}.$$ \rm

\noindent
\bf Remark. \rm This conjecture is now a theorem (joint with George Andrews, in preparation).

\bigskip
\noindent
\bf 13 LOG-CONCAVITY OF CHROMATIC POLYNOMIALS
\bigskip
\noindent
Given a sequence of \bf positive numbers \rm $(a_k)_{k\geq0}$, define the operator $\Cal{L}a_k=a_k^2-a_{k-1}a_{k+1}$. We say $(a_k)_k$ is \it log-concave \rm provided that $\Cal{L}a_k\geq0$ for all $k\geq0$. \it Convenience: \rm enforce that $a_{-1}=0$. 
\smallskip
\noindent
If, after a repeated action of the operator $\Cal{L}$, we find $\Cal{L}^ia_k\geq0$ for $1\leq i\leq m$ and for all $k$, then $(a_k)_k$ is named \it $m$-fold log-concave. \rm The sequence is called \it infinitely log-convex \rm if $\Cal{L}^ia_k\geq0$ for all $i\geq1$.
\bigskip
\noindent
Given a graph $G$ and $x$ distinct colors, denote the number of proper colorings by $\kappa_G(x)$. Then, $\kappa_G(x)$ is referred to as the \it chromatic polynomial \rm of $G$. Here are some examples. 
\bigskip

(N) If $G$ is a graph with $n$ vertices and no edges then $\kappa_G(x)=x^n$.
\bigskip

(T) If $G$ is a tree or a tree with $n$ vertices then $\kappa_G(x)=x(x-1)^{n-1}$.
\bigskip

(K) If $G$ is the complete graph $K_n$ with $n$ vertices then $\kappa_G(x)=x(x-1)\cdots(x-n+1)$.
\bigskip

(C) If $G$ is a cycle with $n$ vertices then $\kappa_G(x)=(x-1)^n+(-1)^n(x-1)$.
\bigskip

(BP) If $G$ is the bipartite graph $K_{n,m}$ then $\kappa_G(x)=\sum_{k=0}^mS(m,k)(x)_k(x-k)^n$, where $S(m,k)$ is 

the \it Stirling number of the second kind \rm and $(x)_k=x(x-1)\cdots(x-k+1)$ is the \it falling factorial. \rm
\bigskip

(CL) If $G$ is the cyclic ladder graph with $2n$ vertices, then 
$$\kappa_G(x)=(x^2-3x+3)^n-(1-x)^{n+1}-(1-x)\,(3-x)^n+(x^2-3x+1).$$

(SB) If $G$ is the signed book graph $B(m,n)$, then
$$\kappa_G(x)=(x-1)^mx^{-n}((x-1)^m+(-1)^m)^n.$$
\bigskip
\noindent
In the present context, June Huh [18] proved the following important result.
\bigskip
\noindent
\bf Theorem. \it The absolute values of coefficients in $\kappa_G(x)$, of any graph $G$, are log-concave. \rm
\bigskip
\noindent
\bf Theoretical. \rm Since polynomials in (N), (T), and (K) are real-rooted, Br\"and\'en's work [16] (conjectured by Fisk [17], McNamara–Sagan [19] and Stanley [20]) shows infinite log-concavity. 
\bigskip
\noindent
\bf Experimental. \rm The coefficients of (C), (BP), (CL), (SB) as well as many such from the long list given by Birkhoff and Lewis [15],  have been tested to be infinitely log-concave.
\bigskip
\noindent
\bf Conjecture 13.1. \it The coefficients of any chromatic polynomial form infinitely log-concave sequences. \rm

\pagebreak

\noindent
\bf 14 DIRECTIVES and UPDATES \rm
\smallskip
\noindent
Enough attempt has been made in recording an accurate account of the results and problems in the
present form of this manuscript. However, the author appreciates any relevant information on the
literature, lemmas or conjectures. These notes will undergo continual changes as more news come;
all such matters (proper documentation, acknowledgment, progress, etc) will be made available at

  $$\text{\tt{http://www.math.tulane.edu/$\sim$tamdeberhan/conjectures.html}}$$

\noindent
Here are a few pointers for some of the conjectures and problems.

For Conjecture 1.8, we are certain that a generating function maneuvering is within reach; for

an interesting combinatorial approach to emulate see [6].

In Conjecture 2.1, the middle two are immediate from [7] that $\sum_nx^n\prod_{u\in\lambda}\frac{t+h_u^2}{h_u^2}=\prod_j(1-x^i)^{-t-1}$;

and the extreme-left and right sides are equal by duality. The main wish here is to prove the rest;

a combinatorial argument is strongly desired.

Problem 2.2 has recently been generalized (to $k$-hooks) and proved in [4].

For Problem 5.3, see [9] for a possible technique on a similar result.

Conjecture 8.1 is amenable to generating function manipulation; yet, a combinatorial 

interpretation is more interesting (ref. [7] is recommended).

Conjecture 11.1 through 11.4 have found resolutions and certain extensions: 

Amdeberhan-Leven \tt{arXiv:1406.2250}\rm; Yang-Zhong-Zhou \tt{arXiv:1406.2583}\rm; 

Xiong \tt{arXiv:1409.7038} \rm and \tt{arXiv:1410.2061}\rm; Nath \tt{arXiv:1411.0339}.  \rm

In regard to Conjecture 12.2 the reader would benefit from consulting one of the most current 

work [14] and the many references therein. This is conjecture is now a theorem (joint G. Andrews).

\enddocument